\documentclass{elsart}
   \usepackage{amsmath}
   \usepackage{subfigure}
   \usepackage[dvips]{graphicx}
   \newcommand{\onehalf}      {\mbox{\small$\frac{1}{2}$}}
   \newcommand{\onethird}     {\mbox{\small$\frac{1}{3}$}}
   \newcommand{\figref}[1]    {\mbox{Fig.~\ref{#1}}}
   \newcommand{\im}{\mbox{$\mathrm{i}$}}
   \renewcommand{\Re}{\mbox{$\mathrm{Re}$}}
   \renewcommand{\Im}{\mbox{$\mathrm{Im}$}}
\journal{Appl. Num. Math.}
\begin{document}
\begin{frontmatter}
\title{Implicit-explicit methods based on\\ strong stability preserving
       multistep time discretizations\thanksref{gsp}}%
\thanks[gsp]{This work was initiated while the author was
           employed by Christian Michelsen Research AS (CMR), Bergen, Norway;
           and supported in part by the sponsors of 
           CMR's Gas Safety Programme (GSP 97--99):
           Agip, British Petroleum, Elf, Exxon, Gaz de France,
           Health and Safety Executive, Mobil, Norsk Hydro,
           Norwegian Petroleum Directorate, Phillips Petroleum, and Statoil.    
      }
\author{Thor Gjesdal}
\address{Norwegian Defence Research Establishment, 
         NO-2027 Kjeller, Norway}
\date{2nd revision, 4 September 2006}
\begin{abstract}
In this note we propose and analyze novel implicit-explicit 
methods based on second order strong stability preserving
multistep time discretizations. 
Several schemes are developed, and a linear stability analysis
is performed to study their properties with respect to the implicit
and explicit eigenvalues. One of the proposed schemes is found
to have very good stability properties, with implicit $A$-stability for
the entire explicit stability domain.
The properties of the other proposed schemes are comparable to those
of traditional methods found in the literature.
\end{abstract}
\begin{keyword}
Implicit-explicit methods, strong stability preserving methods,
advection-diffusion equation, stability
\end{keyword}
\end{frontmatter}
\section{Introduction}
Implicit-explicit schemes (IMEX) are methods for the solution of time-dependent
differential equations in the form
\begin{equation}
\dot{y} = f(y,t) + g(y,t), \label{eq:ode1}
\end{equation}
where $f$ and $g$ are terms with different character such
that one, say $g$, mandates implicit treatment whereas $\dot{y} =f$ can be 
solved efficiently by an explicit method. 
Such systems often arise from spatial (semi-)discretization of time dependent PDEs 
by the method-of-lines
-- a prime example being discretizations of transport equations that may contain
different terms accounting for advection, diffusion, and reaction.

It is beyond the scope of this note to give a complete description of IMEX methods;
for a systematic discussion of multistep IMEX schemes we refer to Ascher et al.~%
\cite{UMAscher_SJRuuth_BTRWetton_1995a}. 
Frank et al.~\cite{JFrank_WHundsdorfer_JGVerwer_1997a} analyzed the stability of 
multistep IMEX schemes, and showed that stable implicit and explicit integrators
do not necessarily lead to a stable implicit-explicit method.
Ascher et al.~\cite{UMAscher_SJRuuth_RJSpiteri_1997a} have also developed IMEX
schemes based on multistage Runge-Kutta integrators.

The Strong Stability Preserving (SSP) property is a generalization of the Total 
Variation Diminishing (TVD) property for hyperbolic problems, and it is a guarantee 
of non-linear stability of the scheme.
Traditional explicit multistep schemes, such as the Adams-Bashforth methods,
do not necessarily have the SSP property~%
\cite{SGottlieb_CWShu_ETadmor_2001a}.
Note however that recent work -- 
by Hundsdorfer and Jaffr\'{e}~\cite{WHundsdorfer_JJaffre_2003a} 
and by Hundsdorfer, Ruuth, and Spiteri~%
\cite{WHundsdorfer_SJRuuth_RJSpiteri_2003a} --
have shown that the traditional multistep methods can be made 
monotonicity-preserving, e.g. with respect to the Total Variation semi-norm,
if a suitable starting procedure is employed.

In this note we discuss the linear stability of  novel IMEX schemes based on  
second order accurate SSP multistep time discretization. 
Although SSP is a non-linear stability property, it is nevertheless appropriate
to study the linear stability of the proposed schemes, since IMEX combinations
do not necessarily inherit the stability properties of the component schemes~%
\cite{JFrank_WHundsdorfer_JGVerwer_1997a}.
A full analysis of the nonlinear stabilities of the proposed schemes is beyond
the scope of this note.

We will show that second order IMEX schemes based on the explicit multistep SSP 
discretizations proposed by Shu~\cite{CWShu_1988a} have stability properties that are
comparable to the Crank-Nicholson/Adams-Bashforth and BDF2 IMEX methods.

\section{Stability analysis}
A general implicit-explicit $k$-step method for \eqref{eq:ode1} can be written in the form
\begin{equation}
\label{eq:general-imex-form}
\sum_{i=0}^{k} a_i y_{n+1-i} = \sum_{i=1}^{k} b_i f_{n+1-i} + \sum_{i=0}^{k} c_i g_{n+1-i}, 
\end{equation}
where $f_n=f(y_n,t_n)$ and $g_n=g(y_n,t_n)$.

The stability analysis below is based on the scalar test equation
\begin{equation}
   \dot{y} = \lambda y + \mu y,
   \label{eq:scalar-test-equation}
\end{equation}
where $\lambda$ and $\mu$ are complex constants that represent the eigenvalues of the 
explicit and implicit operators, respectively.

By applying the multistep method~\eqref{eq:general-imex-form} to the scalar test 
equation~\eqref{eq:scalar-test-equation} and looking for solutions in the
form $y=\zeta^n$, we obtain the characteristic equation
\begin{equation}
\label{eq:general-characteristic-equation}
\sum_{i=0}^{k} a_i \zeta^{n+1-i} 
     - \lambda\sum_{i=1}^{k} b_i \zeta^{n+1-i} - \mu \sum_{i=0}^{k} c_i \zeta^{n+1-i} = 0.
\end{equation}
For convenience, we follow Frank et al.~\cite{JFrank_WHundsdorfer_JGVerwer_1997a} 
in working on the transformed form of the equation 
\begin{equation}
    A(z) - \lambda B(z) - \mu C(z) = 0, \label{eq:final-stability-polynomial}
\end{equation}
where $z=1/\zeta$ and the polynomials $A$, $B$, and $C$ 
correspond to the time derivative, explicit, and implicit operators, respectively,
and are given by
\begin{equation}
      A(z) = \sum_{i=0}^{k} a_i z^i, \quad
      B(z) = \sum_{i=0}^{k} b_i z^i, \quad
      C(z) = \sum_{i=0}^{k} c_i z^i.
      \label{eq:imex-char-polys}
\end{equation}
Note that, since we are working in the transformed variable $z$, the 
scheme~\eqref{eq:general-imex-form} is stable if all the roots of the characteristic 
equation~\eqref{eq:final-stability-polynomial} are in the exterior of the unit disk,
$|z|\geq 1$ (with strict inequality if $z$ is a multiple root).

To investigate the stability properties of the combined IMEX scheme we will
study the image of the exterior of the unit circle under the mapping
\begin{equation}
   \label{eq:phi_lambda}
   \mu = {\varphi}_{\lambda}(z) = \frac{A(z)-\lambda B(z)}{C(z)}.
\end{equation}
The boundary of the stability domain, $\mathcal{S}$,  for the explicit method
can be readily determined  from 
\[ \partial\mathcal{S} = \left\{\lambda : \lambda(\theta)= {A(e^{\im\theta})}/{B(e^{\im\theta})},\quad -\pi\leq\theta<\pi\right\}. \]
Similarly, we can find the stability region, $\mathcal{D}$, of the implicit method
from 
\[ \partial\mathcal{D} = \left\{\mu : \mu(\theta)= {A(e^{\im\theta})}/{C(e^{\im\theta})},\quad -\pi\leq\theta<\pi\right\}. \]
%
\section{IMEX methods based on second order multistep SSP schemes}
\subsection{SSP time discretizations}
The first order accurate forward Euler explicit time integration is strongly
stable for some norm, $\|\cdot\|$, if
\[
\| y_{n} + \Delta t f_{n} \| \leq \| y_{n} \|,
\]
under a suitable time step $\Delta t \leq \Delta t_0$, where $\Delta t_0$
is the largest allowable time step for stability. Higher order
time discretizations are strong stability preserving (SSP) if they
retain this property, possibly for a shorter time step
$\Delta t \leq c\Delta t_0$. The constant $c\leq 1$ is, by convention, called the CFL coefficient~%
\cite{SGottlieb_CWShu_ETadmor_2001a}, and should not be confused with the Courant number
of hyperbolic discretizations.
The SSP time discretizations were originally developed by Shu
\cite{CWShu_1988a}
for the Total Variation norm to maintain positivity and monotonicity in
the solution of hyperbolic equations.

\begin{figure}
\begin{center}
\includegraphics{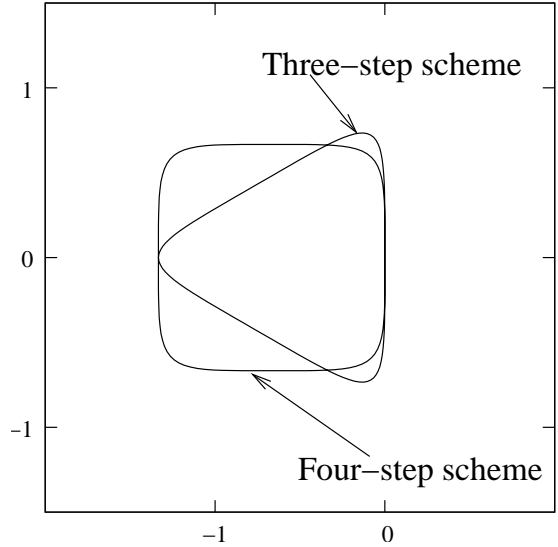}
\end{center}
\caption{\label{fig:explicit-2nd-order-stabreg}
Explicit stability domain, $\mathcal S$, for the two SSP time discretizations
\protect{\eqref{eq:explicit-2nd-order-ssp-schemes}}.
}
\end{figure}

Shu~\cite{CWShu_1988a} proposed a family of explicit SSP multistep schemes.
We will restrict the discussion to second order schemes with positive
coefficients. We thus consider the three- and four-step second order schemes:
\begin{subequations}
\label{eq:explicit-2nd-order-ssp-schemes}
\begin{gather}
\dot{y}_{\Delta t}^{[3]} = \frac{4y_{n+1}-3y_{n}-y_{n-2}}{6\Delta t} = f_{n},
\label{eq:explicit-2nd-order-ssp-schemes-k3}\\
\dot{y}_{\Delta t}^{[4]} = \frac{9y_{n+1}-8y_{n}-y_{n-3}}{12\Delta t} = f_{n},
\label{eq:explicit-2nd-order-ssp-schemes-k4}
\end{gather}
where $\dot{y}_{\Delta t}^{[k]}$ denotes the $k$-step discretized time derivative
for the SSP scheme.
\end{subequations}
These methods are strong stability preserving with  CFL coefficients 
$c=1/2$ and $c=2/3$ for $k=3$ and $k=4$, respectively.
The polynomials~\eqref{eq:imex-char-polys} that characterize these schemes
are
\begin{subequations}
\label{eq:A-poly}
\begin{gather}
A^{[3]}(z) = (4-3z-z^3)/6, \label{eq:A-poly-a}\\
A^{[4]}(z) = (9-8z-z^4)/12, \label{eq:A-poly-b}
\end{gather}
\end{subequations}
and 
\begin{equation}
B(z) = z. \label{eq:B-poly}
\end{equation}
We show the explicit stability domain, $\mathcal{S}$, for the two schemes~%
\eqref{eq:explicit-2nd-order-ssp-schemes}
in \figref{fig:explicit-2nd-order-stabreg}.
Notice that both schemes have an appreciable part of the stability region close 
to the imaginary axis, something that is advantageous for the solution of 
hyperbolic terms.
\subsection{Construction of implicit integrators}
\label{sec:implicit}
To develop IMEX methods based on the multistep SSP schemes, we 
first construct implicit schemes by replacing the term  $f(y^n)$ 
in \eqref{eq:explicit-2nd-order-ssp-schemes} 
above by linear combinations of the implicit operator evaluated at the available
time levels
\[
   \sum_{j=0}^{k} c_j g_{n+1-j},
\]
such that the resulting formula for calculating $y^{n+1}$ is second order accurate at time level $n$.
We must obviously require $c_0\neq 0$ for the scheme to be implicit, and furthermore
that $\sum c_j = 1$.

Evaluating the order conditions we obtain an under-determined system that allows
infinitely many solutions with several free parameters. We will not go into all the
possibilities here. Instead we focus on the straightforward approximations
\begin{subequations}
\label{eq:implicit-integrators}
\begin{gather}
\dot{y}_{\Delta t}^{[k]} =  \onethird \bigg( 2g_{n+1} + g_{n-2} \bigg),
       \label{eq:implicit-integrators-biased} \\
\dot{y}_{\Delta t}^{[k]} =  
      \onehalf \bigg( (1-\beta)g_{n+1} + 2\beta g_{n} + (1-\beta)g_{n-1} \bigg), 
       \label{eq:implicit-integrators-centred} 
\end{gather}
\end{subequations}
where $\beta<1$ is a non-negative algorithmic parameter.

The polynomials~\eqref{eq:imex-char-polys} that correspond to the right-hand-sides
in \eqref{eq:implicit-integrators} are then
\begin{subequations}
\label{eq:C-poly}
\begin{gather}
C_1(z) = (2+z^3)/3, \label{eq:C-poly-a}\\
C_2(z) = \left((1-\beta) +2\beta z + (1-\beta) z^2 \right)/2. \label{eq:C-poly-b}
\end{gather}
\end{subequations}

An implicit method has $A$-stability if it is stable for all eigenvalues, $\mu$, in the
right half plane, $\Re(\mu)\leq 0$, and it has $A(\alpha)$-stability if it is stable
for eigenvalues in a wedge-shaped domain about the negative real axis 
with half-angle $\alpha$.

\begin{figure}
\begin{center}
\includegraphics{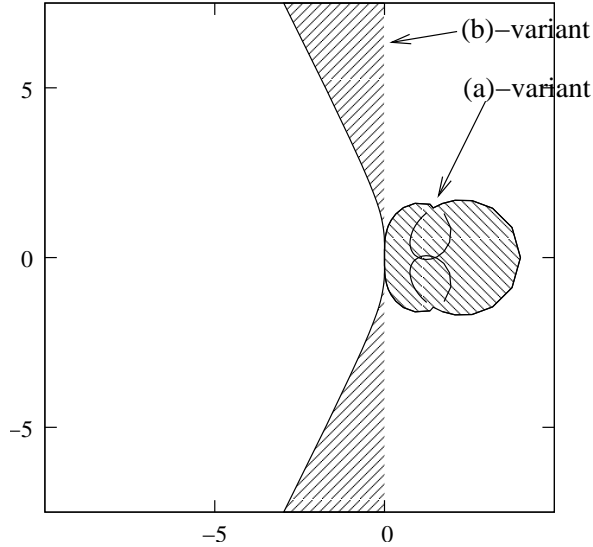}
\end{center}
\caption{\label{fig:implicit-integrators-stabreg}
Stability regions, $\mathcal{D}$, for the implicit integrators 
\protect{\eqref{eq:implicit-integrators}}, where $k=3$ and $\beta=0$.
The methods are stable outside the shaded regions.
}
\end{figure}

In \figref{fig:implicit-integrators-stabreg} we show the stability regions of 
the $k=3$ versions of the implicit integrators~%
\eqref{eq:implicit-integrators}.
Note that the scheme~\eqref{eq:implicit-integrators-biased} appears to be $A$-stable, 
whereas the other variant~\eqref{eq:implicit-integrators-centred} is $A(\alpha)$-stable. 
We can show that these stability properties carry over to the $k=4$ versions of the schemes. 
That is, \eqref{eq:implicit-integrators-biased} is $A$-stable
and~\eqref{eq:implicit-integrators-centred} is $A(\alpha)$-stable, for both
$k=3$ and $k=4$.  

\begin{lem}
\label{lemma:m3p2-implicit-b-stability}
The implicit method~\eqref{eq:implicit-integrators-biased} is $A$-stable for both
$k=3$ and $k=4$.
\end{lem}
\begin{pf}
Obviously, for a purely implicit method $\lambda=0$,
and we consider the mapping
\[
\varphi_0(z)={A^{[k]}(z)}/{C_1(z)}.
\]
$A$-stability requires that the image of the unit circle under this mapping 
is entirely in the right half plane, i.e.~$\Re\left(\varphi_0(e^{\im\theta})\right)\geq 0$.

(i) For $k=3$ we find
\begin{gather}
\Re\left(\varphi_0(e^{\im\theta})\right) = \frac{14 - 12\cos\theta - 6\cos2\theta + 4\cos3\theta}
                                    {4\left|2+e^{\im3\theta}\right|^2} \\
                             = \frac{5 - 6\cos\theta - 3\cos^2\theta + 4\cos^3\theta}
                                    {\left|2+e^{\im3\theta}\right|^2}.
\end{gather}
The denominator in this expression is obviously positive, we must therefore
show that the numerator is non-negative. Let $x=\cos\theta$ and write
\[ f(x) = 4x^3 - 3x^2 -6x + 5 = (4x-5)(x-1)^2. \]
It is straightforward to verify that $f(x)$ is non-negative for $-1\leq x \leq 1$, 
and we have shown that \eqref{eq:implicit-integrators-biased} is $A$-stable for 
$k=3$.

(ii) Similarly, for $k=4$, we have
\[
\varphi_0(e^{\im\theta})={A^{[4]}(e^{\im\theta})}/{C_1(e^{\im\theta})},
\]
leading to
\[
\Re\left(\varphi_0(e^{\im\theta})\right) = \frac{6 - 11\cos\theta + 9\cos^3\theta - 4\cos^4\theta}{4\left|2+e^{\im3\theta}\right|^2}.
\]
Again we let $x=\cos\theta$ and write the numerator as
\[f(x) = -4x^4-9x^3-x+6 = -(4x^2-x-6)(x-1)^2, \]
which we readily verify is non-negative for $-1\leq x \leq 1$.
\qed
\end{pf}
\begin{lem}
\label{lemma:m3p2-implicit-c-stability} The implicit method ~\eqref{eq:implicit-integrators-centred} is $A(\alpha)$-stable,
with $\alpha=\alpha_0$ given by
\begin{subequations}
\label{eq:implicit-centred-tan-alpha-0}
\begin{gather}
\tan \alpha_0 = \frac{(2+\gamma)\sqrt{1-\gamma^2}}{(\gamma-1)^2}, \quad\mbox{for } k=3, 
\label{eq:implicit-centred-tan-alpha-0-k3}\\
\tan \alpha_0 = \frac{(2+\gamma^2)\sqrt{1-\gamma^2}}{2-3\gamma+\gamma^3}, \quad\mbox{for } k=4,
\label{eq:implicit-centred-tan-alpha-0-k4} 
\end{gather}
\end{subequations}
where $\gamma=\beta/(\beta-1)$.
\end{lem}
\begin{pf}
Again, we let $\lambda=0$,
and we consider the mapping
\[
\varphi_0(e^{\im\theta})={A^{[k]}(e^{\im\theta})}/{C_2(e^{\im\theta})},
\]
shown in \figref{fig:implicit-integrators-stabreg}, 
to determine the half-angle, $\alpha_0$.
Note that $\varphi_0(\exp(\im\theta))$ is singular in the limit 
$\cos\theta\rightarrow \gamma=\beta/(\beta-1)$. 
We only consider $\beta\leq 1/2$; since for $\beta>1/2$ the scheme behaves more like an explicit method.
The angle $\alpha_0$ is given by the slope of the asymptote, which can be written
\[
\tan\alpha_0 
             = \lim_{\cos\theta\rightarrow\gamma}\frac{\Im(\varphi_0(e^{\im\theta}))}{\Re(\varphi_0(e^{\im\theta}))},
\]
We then have
\begin{gather*}
\tan\alpha_0 
   = \frac{-4\sin\theta^*-\sin 3\theta^*}{4-3\cos\theta^*-\cos 3\theta^*}, \quad\mbox{for } k=3, \\
\tan \alpha_0  
   = \frac{-8\sin\theta^*-\sin 4\theta^*}{9-8\cos\theta^*-\cos 4\theta^*}, \quad\mbox{for } k=4, 
\end{gather*}
from which the  results follow.
\qed
\end{pf}

The natural choice $\beta=0$, i.e.~$\gamma=0$, gives
$\tan\alpha_0=2$ and $\tan\alpha_0=1$ for the three-step and the four-step schemes
respectively. In both cases, this is the optimal parameter value that maximizes $\alpha_0$.

\subsection{Stability of the IMEX schemes}
We can the construct implicit explicit methods by combining the above
SSP explicit schemes~%
\eqref{eq:explicit-2nd-order-ssp-schemes} 
and the implicit integrators~
\eqref{eq:implicit-integrators}.
Thus we consider 
\begin{subequations}
\label{eq:imex-schemes}
\begin{gather}
\dot{y}_{\Delta t}^{[k]} =  f(y^{n}) +\onethird \bigg( 2g(y^{n+1}) + g(y^{n-2}) \bigg).
      \label{eq:imex-scheme-biased} \\
\dot{y}_{\Delta t}^{[k]} =  f(y^{n}) +
      \onehalf \bigg( (1-\beta)g(y^{n+1}) + 2\beta g(y^{n}) + (1-\beta)g(y^{n-1}) \bigg), 
      \label{eq:imex-scheme-centred} 
\end{gather}
\end{subequations}

Frank et al.~\cite{JFrank_WHundsdorfer_JGVerwer_1997a} presented two general
stability results for multistep IMEX schemes. The first gives restrictions
on the explicit eigenvalues in order to retain the full stability domain of
the implicit operator, whereas the second gives restrictions on the implicit 
eigenvalues in order to retain the full stability domain of the explicit operator.
Unfortunately, neither of these results appear to be applicable to the proposed 
schemes, and we will therefore study the stability properties
of~\eqref{eq:imex-schemes} by working directly on the mapping 
$\varphi_{\lambda}(\exp(\im\theta))$ from \eqref{eq:phi_lambda}.
Remember that the IMEX schemes are stable with respect to implicit eigenvalues, $\mu$,
that are in the exterior of the mapping of the unit disk under $\varphi_{\lambda}$.
For a method to be $A$-stable, the entire image of the unit circle must be in the right
half-plane.

\subsubsection{Stability of scheme~\protect{\eqref{eq:imex-scheme-biased}}}
In this section we will discuss the stability of 
the scheme~\eqref{eq:imex-scheme-biased}, characterized by the
polynomials~\eqref{eq:A-poly}, \eqref{eq:B-poly}, and \eqref{eq:C-poly-a}.
For the $k=3$ variant of the scheme we propose the following stability result:
\begin{conj}
The scheme~\eqref{eq:imex-scheme-biased} with $k=3$ is $A$-stable with respect to $\mu$ 
for $\lambda\in\mathcal{S}$.
\end{conj}
Recall that $A$-stability requires that the image of the unit circle under the
mapping $\varphi_{\lambda}(z)$ must be entirely in the right half-plane for 
$\lambda\in\mathcal{S}$, i.e.
\[ \Phi(\lambda,z) = \Re\left(\varphi_{\lambda}(e^{\im\theta})\right) \geq 0, \quad-\pi\leq\theta < \pi, \quad\lambda\in\mathcal{S}. \]
As a consequence of the maximum principle, $\Phi$ will attain its minimal value
w.r.t.~$\lambda$ on the boundary of the explicit stability domain.
We need therefore only consider $\lambda\in\partial\mathcal{S}$.

\begin{figure}
\begin{center}
{\includegraphics[scale=1.0]{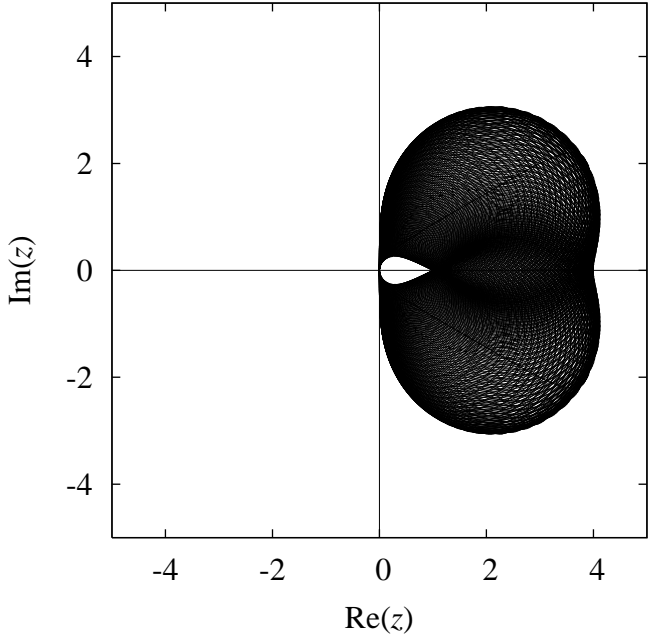}}
\end{center}
\caption{Image of the unit circle under the mapping $\varphi_{\lambda}(z)$, 
with $\lambda \in \partial\mathcal{S}$,
for the schemes given by~\eqref{eq:imex-scheme-biased} with $k=3$. 
The implicit stability domain of the method is outside the shaded region.
\label{fig:imex_bias_philambda3}
}
\end{figure}

To support the assertion that \eqref{eq:imex-scheme-biased} with $k=3$ is $A$-stable, 
we first show, in~\figref{fig:imex_bias_philambda3}, the image of the unit circle under 
$\varphi_{\lambda}$ for $\lambda \in \partial \mathcal{S}$. 
Furthermore, let $\lambda^*=A^{[3]}(\exp(\im\theta^*))/B(\exp(\im\theta^*))$ be an 
eigenvalue on the boundary of the explicit stability domain.
Then $\varphi_{\lambda^*}(\exp(\im\theta^*))=0$, and we can investigate the
behaviour of $\varphi_{\lambda^*}$ close to this zero by performing a Taylor
expansion.
To leading order we obtain
\begin{equation}
\Re(\varphi_{\lambda^*}(e^{\im\theta})) 
\approx \frac{1-\cos 3\theta^*}{5+4\cos 3\theta^*}(\theta-\theta^*)^2 \geq 0.
\end{equation}
We have neither analytically nor numerically managed to find other zeros or 
negative values for $\Phi(\lambda,z)$ and
the scheme does therefore indeed appear to be $A$-stable.

\begin{figure}
\begin{center}
{\includegraphics[scale=1.0]{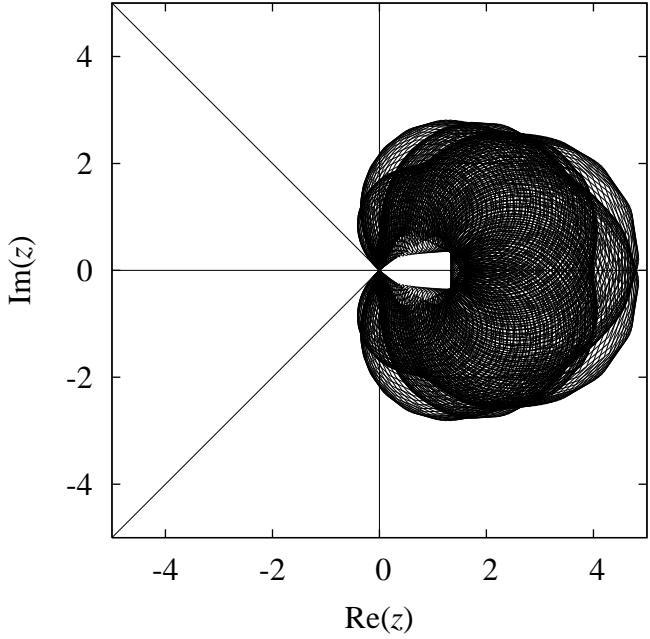}}
\end{center}
\caption{Image of the unit circle under the mapping $\varphi_{\lambda}(z)$, 
with $\lambda \in \partial\mathcal{S}$,
for the schemes given by~\eqref{eq:imex-scheme-biased} with $k=4$. 
The implicit stability domain of the method is outside the shaded region.
\label{fig:imex_bias_philambda4}
}
\end{figure}

For the $k=4$ variant we propose:
\begin{conj}
The scheme~\eqref{eq:imex-scheme-biased} with $k=4$ is $A(\alpha)$-stable 
with respect to $\mu$ for $\lambda\in\mathcal{S}$. 
The half-angle, $\alpha$, of the wedge of stability is $\alpha\approx 0.23\pi$.
\end{conj}
In~\figref{fig:imex_bias_philambda4} we show the image of the unit circle under 
$\varphi_{\lambda}$ for $\lambda \in \partial \mathcal{S}$, and we observe that
the method version appears to be $A(\alpha)$ stable with $\alpha\approx\pi/4$.
Carrying out a Taylor expansion of $\varphi_{\lambda^*}$ as above, we find to 
leading order
\begin{subequations}
\begin{gather}
\Re(\varphi_{\lambda}(e^{\im\theta})) \approx \frac{3}{4}\;\frac{\sin\theta^* - 3\sin 3\theta^* + 2\sin 4\theta^*}
                                                    {\sin^2 3\theta^* + ( \cos 3\theta^*+2)^2}(\theta-\theta^*),          \\
\Im(\varphi_{\lambda}(e^{\im\theta})) \approx \frac{3}{4}\;\frac{6+\cos \theta^* + 3\cos 3\theta^* + 2\cos 4\theta^*}
                                                    {\sin^2 3\theta^* + ( \cos 3\theta^*+2)^2}(\theta-\theta^*).
\end{gather}
\end{subequations}
We can estimate the angle for $A(\alpha)$-stability from the slope of the first order approximation:
\begin{equation}
\tan \alpha \approx \inf_{\theta^*} \left| 
     \frac{6+\cos \theta^* + 3\cos 3\theta^* + 2\cos 4\theta^*}%
          {\sin\theta^* - 3\sin 3\theta^* + 2\sin 4\theta^*}
     \right|.
\end{equation}
Numerically we find that $\tan \alpha \approx 0.89$, and this gives $\alpha\approx 0.23\pi$ which corresponds well with
the stability domain shown in \figref{fig:imex_bias_philambda4}.
\subsubsection{Stability of scheme~\protect{\eqref{eq:imex-scheme-centred}}}
We then consider the scheme~\eqref{eq:imex-scheme-centred}, characterized by the
polynomials~\eqref{eq:A-poly}, \eqref{eq:B-poly}, and \eqref{eq:C-poly-b}.
For this scheme, we can show the following stability result:
\begin{lem}
\label{lemma:m3p2-imex-stability}
The implicit-explicit method~\eqref{eq:imex-scheme-centred} 
is $A(\alpha)$-stable with respect to the implicit eigenvalues, $\mu$, 
for $\beta\leq\onehalf$ if
  \[
   \lambda\in \mathcal{S}^{\nu}= \mathcal{S}\cap\{z:|\Im(z)|\leq \nu \},\quad \mbox{and}\quad
   |\alpha|\leq\alpha_{\nu} 
  \]
where $\nu$ and $\alpha_{\nu}$ are given by 
\begin{subequations}
\label{eq:tana-nu-k=3}
\begin{gather}
\nu \leq (2+\gamma)\sqrt{1-\gamma^2}/3, \label{eq:lambda-i-bound-k3}\\
\tan\alpha_{\nu} =
    \frac{(2+\gamma)\sqrt{1-\gamma^2} - 3\nu }{(\gamma-1)^2}, \label{eq:tan-a-bound-k3}
\end{gather}
\end{subequations}
for $k=3$, and 
\begin{subequations}
\label{eq:tana-nu-k=4}
\begin{gather}
\nu \leq (2+\gamma^2)\sqrt{1-\gamma^2}/3, \label{eq:lambda-i-bound-k4}\\
\tan \alpha_{\nu} = \frac{(2+\gamma^2)\sqrt{1-\gamma^2} -3\nu}{2-3\gamma+\gamma^3}.
\label{eq:tan-a-bound-k4}
\end{gather}
\end{subequations}
for $k=4$.
\end{lem}
\begin{pf}
Recall from Sec.~\ref{sec:implicit} that, in this case, the mapping $\varphi_{\lambda}(\exp(\im\theta))$ is singular in
the limit $\cos\theta\rightarrow \gamma=\beta/(\beta-1)$. 
Taking into account a non-zero $\lambda=\lambda_r+\im\lambda_i\in\mathcal{S}$ and proceeding as in the analysis of
Sec.~\ref{sec:implicit}, we have
\begin{gather*}
\tan\alpha_{\lambda} = \inf_{\lambda}\left( 
   \frac{(2+\gamma)\sqrt{1-\gamma^2}-3\lambda_i}{(\gamma-1)^2 + 3\lambda_r}\right), \qquad \mbox{for $k=3$} \\
\tan\alpha_{\lambda} = \inf_{\lambda}\left( 
   \frac{(2+\gamma^2)\sqrt{1-\gamma^2}-3\lambda_i}{\gamma^3 -3\gamma + 2 + 3\lambda_r}\right) \qquad \mbox{for $k=4$}.
\end{gather*}
We observe that for large $|\lambda_i|$, say $|\lambda_i|>\nu$; $\tan\alpha_{\lambda}$ may become
negative. This means that the wedge of stability will be rotated entirely into the upper or lower left quadrants. 
In order to have a useful stability domain that -- at least -- contains the negative real axis, 
$\tan\alpha_{\lambda}$ must obviously be non-negative, and
we must introduce an upper limit on the imaginary part of the explicit eigenvalues, 
\begin{gather*}
|\lambda_i|\leq \nu \leq (2+\gamma)\sqrt{1-\gamma^2}/3, 
\qquad\mbox{for $k=3$}\\
|\lambda_i|\leq \nu \leq (2+\gamma^2)\sqrt{1-\gamma^2}/3, 
\qquad\mbox{for $k=4$}.
\end{gather*}
We obtain a lower bound on the angle $\alpha$ if
we choose $\lambda = \pm \im\nu$, and we have:
\begin{gather*}
\tan\alpha_{\nu} =
    \frac{(2+\gamma)\sqrt{1-\gamma^2} - 3\nu }{(\gamma-1)^2}, 
\\
\tan \alpha_{\nu} = \frac{(2+\gamma^2)\sqrt{1-\gamma^2} -3\nu}{2-3\gamma+\gamma^3}.
\end{gather*}
for $k=3$ and $k=4$, respectively.%
\qed
\end{pf}

\begin{figure}
   \begin{center}
      \includegraphics[]{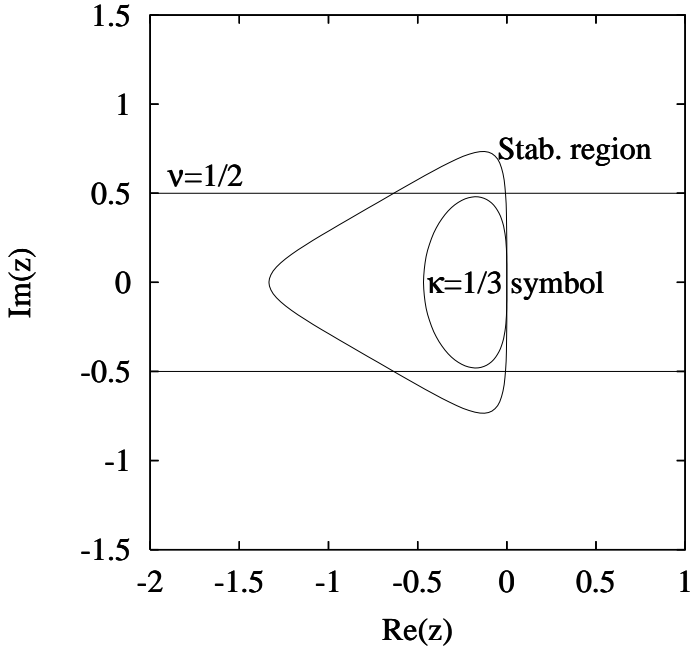} 
   \end{center}
   \caption{Explicit stability region $\mathcal{S}$ for the time discretization~%
\protect{\eqref{eq:explicit-2nd-order-ssp-schemes-k3}}, with the Fourier symbol of the third order 
            $\kappa=1/3$ finite difference advection scheme ($\sigma=0.35$), 
            and the bound ($\gamma=0.5$) that defines the restricted stability 
            region 
            $\mathcal{S}^{\gamma}$ in Lemma~\ref{lemma:m3p2-imex-stability}.
   \label{fig:m3p2-gamma-kappa}
   }
\end{figure}

Note that the natural choice $\beta=0$ maximizes both the angle $\alpha$ and the 
upper bound on $\nu$, and thus gives the largest stability domains
for both the implicit and explicit operators.
Furthermore, although the above results require that we use  a reduced explicit 
stability domain, the restrictions are not prohibitive. 
The choice of the bound $\nu$ can be guided by the discretization of the explicit
operator.
As an example, we compare 
-- in~\figref{fig:m3p2-gamma-kappa} --
the Fourier symbol of the third order finite difference
$\kappa=1/3$ advection scheme 
(see for example~\cite{hundsdorfer-koren-etal:advection-jcp-95})
with the reduced stability domain,
$\mathcal{S}^{0.5}$  for $k=3$.
The Courant number in this particular example is $\sigma=0.35$ which is well above the
value for which the scheme is expected to be SSP ($\sigma_{\mathrm{SSP}}=0.25$).
In the table below, we will use the value $\nu=1/3$ which corresponds fairly well to
$\sigma=0.25$ for this advection discretization.

\section{Discussion and concluding remarks}
In this note we have proposed and analyzed the linear stability properties of implicit-explicit
methods based on the Strong Stability Preserving multistep methods introduced by Shu~%
\cite{CWShu_1988a}.
Of the four variant methods that we consider, one scheme stands out from the rest.
The three-step scheme~\eqref{eq:imex-scheme-biased} with $k=3$ retains full implicit 
$A$-stability for the entire explicit stability domain. 
In fact, 
we are not aware of any other second order IMEX scheme with this property.
For the other methods, we have 
shown $A(\alpha)$-stability for the entire, or slightly reduced explicit stability domains.

It is instructive to compare with the traditional multistep IMEX methods such as the
modified Crank-Nicholson/Adams-Bashforth (mCNAB) scheme and the IMEX BDF2 scheme.
The mCNAB scheme introduces an algorithmic parameter, $c$, to the implicit integrator 
and can be written
\[
y_{n+1}-y_{n} = \onehalf \Delta t\big( 3f_{n} - f_{n-1} + (1+c)g_{n+1} + (1-2c)g_{n} + cg_{n-1} \big).
\]
Note that $c=0$ recovers the original Crank-Nicholson scheme for the implicit integrator, 
whereas other choices of $c$ may improve the stability properties of the method, 
see~\cite{UMAscher_SJRuuth_BTRWetton_1995a,JFrank_WHundsdorfer_JGVerwer_1997a} 
for a discussion.
The IMEX BDF2 scheme is given by
\[
3y_{n+1} -4y_{n} +y_{n-1} = 2\Delta t\big(2f_{n}-f_{n-1} + g_{n+1}   \big).
\]
Both these schemes are restricted to implicit $A(\alpha)$-stability if we require
full explicit stability~%
\cite{JFrank_WHundsdorfer_JGVerwer_1997a}.
We summarize the implicit stability properties for all these methods in 
Table~\ref{table:alpha-table}.
Furthermore, all the methods appear to have comparable non-linear stability properties. Whereas the schemes
proposed in this work are based on strong stability preserving methods, both the Adams-Bashforth method 
and the extrapolated BDF2 method was recently shown to have good monotonicity-preserving properties provided
suitable starting procedures were used~%
\cite{WHundsdorfer_JJaffre_2003a,WHundsdorfer_SJRuuth_RJSpiteri_2003a}.
Finally, it is difficult to assess the computational cost of each method, as this may depend strongly on details of
the implementation. It is however reasonable to assume that all the methods considered in this study as well
as the mCNAB and IMEX BDF2 schemes have comparable computational cost. The memory requirements may however
be larger for the SSP-based schemes as data from more time levels need to be stored.
\begin{table}
\begin{center}
\caption{Comparison of the implicit stability domains with full explicit stability
         for the schemes considered in this work and the traditional IMEX BDF2 and mCNAB methods.
         The results for the latter methods are taken from Reference~%
\protect{\cite{JFrank_WHundsdorfer_JGVerwer_1997a}}.
        \label{table:alpha-table}
}
\begin{tabular}{lc}
Scheme & $\alpha$ \\ \hline
\protect{\eqref{eq:imex-scheme-biased}} ($k=3$)                       & $\pi/2$   \\
\protect{\eqref{eq:imex-scheme-biased}} ($k=4$)                       & $0.23\pi$ \\
\protect{\eqref{eq:imex-scheme-centred}} ($\beta=0$, $\nu=1/3$, $k=3$) & $\pi/4$    \\
\protect{\eqref{eq:imex-scheme-centred}} ($\beta=0$, $\nu=1/3$, $k=4$) & $0.15\pi$ \\
IMEX BDF2                                                             & $0.31\pi$   \\
mCNAB ($c=0$)                                                          & $0$         \\
mCNAB ($c=1/8$)                                                        & $0.12\pi$   \\
mCNAB ($c=1/2$)                                                        & $0.23\pi$   \\ \hline
\end{tabular}
\end{center}
\end{table}

We will also briefly discuss IMEX Runge-Kutta methods. 
Ascher et al.~\cite{UMAscher_SJRuuth_RJSpiteri_1997a}
suggested that a two-stage, second-order IMEX Runge-Kutta scheme based on the trapezoidal rule
could be an alternative to the CNAB method in  `most situations'. 
There are also schemes in this family with very good SSP proprieties.
Analysis of IMEX RK2 schemes shows unfortunately that there are severe restrictions on the explicit 
stability domain, in particular for implicit eigenvalues with large magnitude~%
\cite{WHundsdorfer_JGVerweer_2003a}.
Higher order IMEX RK methods have been proposed and analyzed by
Ascher et al.~\cite{UMAscher_SJRuuth_RJSpiteri_1997a},
Pareschi and Russo~\cite{LPareschi_GRusso_2000a},
and Higueras~\cite{IHigueras_2004a}, but these have not been considered in the present work.

\bibliographystyle{plain}
\bibliography{references}

\begin{thebibliography}{10}

\bibitem{UMAscher_SJRuuth_RJSpiteri_1997a}
U.~M. Ascher, S.~J. Ruuth, and R.~J. Spiteri.
\newblock Implicit-explicit {Runge-Kutta} methods for time-dependent partial
  differential equations.
\newblock {\em Appl. Numer. Math.}, 25:151--167, 1997.

\bibitem{UMAscher_SJRuuth_BTRWetton_1995a}
U.~M. Ascher, S.~J. Ruuth, and B.~T.~R. Wetton.
\newblock Implicit-explicit methods for time-dependent partial differential
  equations.
\newblock {\em SIAM J. Numer. Anal.}, 32:797--823, 1995.

\bibitem{JFrank_WHundsdorfer_JGVerwer_1997a}
J.~Frank, W.~Hundsdorfer, and J.~G. Verwer.
\newblock On the stability of implicit-explicit linear multistep methods.
\newblock {\em Appl. Numer. Math.}, 25:193--205, 1997.

\bibitem{SGottlieb_CWShu_ETadmor_2001a}
S.~Gottlieb, C.~W. Shu, and E.~Tadmor.
\newblock Strong stability-preserving high-order time discretization methods.
\newblock {\em SIAM Rev.}, 43:89--112, 2001.

\bibitem{IHigueras_2004a}
I.~Higueras.
\newblock Strong stability for additive {Runge-Kutta} methods.
\newblock In {\em Book of abstracts. {ICOSAHOM 2004}}. Brown University, 2004.

\bibitem{WHundsdorfer_JJaffre_2003a}
W.~Hundsdorfer and J.~Jaffr\'{e}.
\newblock Implicit-explicit time stepping with spatial discontinuous finite
  elements.
\newblock {\em Appl. Numer. Math.}, 45:231--154, 2003.

\bibitem{hundsdorfer-koren-etal:advection-jcp-95}
W.~Hundsdorfer, B.~Koren, M.~van Loon, and J.~G. Verwer.
\newblock A positive finite-difference advection scheme.
\newblock {\em J. Comput. Phys.}, 117:35--46, 1995.

\bibitem{WHundsdorfer_SJRuuth_RJSpiteri_2003a}
W.~Hundsdorfer, S.~J. Ruuth, and R.~J. Spiteri.
\newblock Monotonicity-preserving linear multistep methods.
\newblock {\em SIAM J. Num. Anal.}, 41:605--623, 2003.

\bibitem{WHundsdorfer_JGVerweer_2003a}
W.~Hundsdorfer and J.~G. Verweer.
\newblock {\em Numerical Solution of Time-Dependent
  Advection-Diffusion-Reaction Equations}.
\newblock Springer, 2003.

\bibitem{LPareschi_GRusso_2000a}
L.~Pareschi and G.~Russo.
\newblock Implicit-explicit {Runge-Kutta} schemes for stiff systems of
  differential equations.
\newblock In L.~Brugnano and D.~Trigiante, editors, {\em Recent Trends in
  Numerical Analysis}, volume~3 of {\em Advances in Computation: Theory and
  Practice}, pages 269--289. Nova Science, 2000.

\bibitem{CWShu_1988a}
C.~W. Shu.
\newblock Total-variation-diminishing time discretizations.
\newblock {\em SIAM J. Sci. Stat. Comp.}, 9:1073--1084, 1988.

\end{thebibliography}
\end{document}